# REDUCING THE AXIOMS OF HYPERGROUPS, HYPERFIELDS, HYPERMODULES AND RELATED STRUCTURES: A NEW AXIOMATIC BASIS FOR HYPERCOMPOSITIONAL STRUCTURES

CHRISTOS G. MASSOUROS

ABSTRACT: This paper is concerned with the axiomatic basis of structures within Hypercompositional Algebra. It is proven that the axioms employed in the definition of numerous hypercompositional structures lack independence. Accordingly, novel definitions are introduced in this work which minimize the established definitions by reducing the necessary set of axioms.

## 1. Introduction

The independence of axioms, and the question of whether some are derivable from others, lies at the heart of mathematical thought. Since antiquity, the independence of the parallel postulate from the remaining axioms of Euclidean geometry has been a source of profound reflection. This centuries–long quest eventually gave rise to entirely new worlds of thought: the hyperbolic geometry of Lobachevsky (1829) and the elliptic geometry of Riemann (1854). Such historical developments remind us that the precise formulation of axiomatic systems is not merely a matter of technical detail, but a foundation upon which the scope and shape of whole mathematical theories may depend.

This paper examines the independence of axioms in the context of hypergroups, hyperfields, hypermodules and related structures. It is shown that certain axioms commonly included in their definitions are in fact derivable from the others and hence need not be postulated independently. This observation not only simplifies and clarifies the foundations of hypercompositional structures but also maintains a strict conceptual boundary between postulates and theorems.

Section 2 begins with the study of hypergroups, the foundational hypercompositional structure introduced in 1934 by the French mathematician F. Marty [1]. Despite his all-too-brief life—as Marty was killed in action during the Second World War—he published three seminal papers on the subject [1–3]. These works laid the groundwork for hypergroups and the subsequent development of more complex hypercompositional structures. The excerpt from Marty's paper [1] in which the definition of a hypergroup is given, is reproduced verbatim below:

*Soit un ensemble d'éléments, non vide, possédant quatre lois de combinaison $AB$, $BA$, $\dfrac{A}{B}$, $\dfrac{A}{|B}$;*

*chacune d'elles pouvant avoir plusieurs déterminations; les deux premières étant associatives. Si C est une détermination de AB, nous écrirons  AB ⊃ C  (AB contient C). Nous dirons que la famille constitue un hypergroupe si les deux divisions sont liées à la multiplication par les relations suivantes:*

$$\frac{A}{B} \supset C \ \rightleftarrows \ BC \supset A$$

$$et \qquad \frac{A}{|B} \supset C \ \rightleftarrows \ CB \supset A$$

One year later, in 1935, F. Marty, in [2], studied a special class of hypergroups possesing a neutral element, which he termed "*normal hypergroups*" (*hypergroupes normaux*). In 1936, he repeated the above definition of a hypergroup in [3]. However, regardless of which book or article one consults today, this is not the definition of a hypergroup that appears in the literature. The formulation of the definition encountered in all textbooks and articles on hypergroups is due to M. Krasner. Yet this fact is nowhere explicitly mentioned; on the contrary, the definition is erroneously attributed to F. Marty. A careful study of the literature reveals that in 1937, in his work [4], Krasner provided the following definition of the hypergroup:

*DÉFINITION 1: Un ensemble K organise par une loi de composition ab de chaque couple a, b de ses éléments s'appelle hypergroupe par rapport à cette loi de composition si*
$1^0$    *ab est un sous-ensemble non-vide de K.*
$2^0$    *(ab)e = a(bc) (loi associative).*
$3^0$    *Pour chaque couple a,c ∈ H il existe x ∈ H tel que c ∈ ax et il existe x′ ∈ H tel que c ∈ x′a*

Two years later, in 1939, Krasner published in [5] the definition of the hypergroup in the form that is used today:

*Un ensemble H, ou est définie une loi de composition de ses éléments, s'appelle hypergroupe par rapport a cette loi si les trois conditions suivantes sont vérifiées:*
1.    *Le compose ab d'un b∈H par un a∈H est un sous-ensemble jamais vide de H.*
2.    *(ab)c = a(bc) (a, b, c ∈H) (loi d'associativité).*
3.    *Pour tout a∈H, aH =Ha= H.*

It is worth noting that the above definitions of a hypergroup include three axioms. However, the first axiom — namely, the requirement that the product of any two elements be a non-empty subset — no longer appears in the definition of the hypergroup since it has been incorporated into the very notion of hypercomposition. Indeed, over a non-empty set *H*, a *hypercomposition* (or *hyperoperation*) is defined as a mapping from the Cartesian product *H×H* into the power set of *H*, excluding the empty set (see, e.g. [6-10]).

In Section 2, it is shown that in hypergroups Axiom 1 is redundant, as it follows from Axioms 2 and 3. Consequently, there is no need to impose the usual restriction on the codomain in the definition of the hypercomposition. The same reasoning applies to analogous structures such as right and left almost hypergroups [11] and H$_v$-groups [12], as well as to more complex hypercompositional structures built upon them, such as hyperrings [13], hyperfields [13], hyperringoids [14], hypermodules [15], hypermoduloids [16], and others [17–22]. Section 3 focuses on a particular hypergroup, namely the polysymmetrical hypergroup [22], and demonstrates that the axiom of reversibility can be derived from the remaining axioms.

Section 4 establishes that the same result holds for hyperfields: the axiom of reversibility is not independent but follows directly from the other defining axioms. This result applies equally to unitary hyperrings. Moreover, this dependence is not confined to these structures; it extends naturally to hypermodules and vector hyperspaces as well.

## 2. Reduction of Axioms in Hypergroups and Related Hypercompositional Structures

In *Éléments de Mathématique, Algèbre* [23], Nicolas Bourbaki introduced the term *magma* (from the Greek verb μάσσω, meaning to knead) to denote a set equipped with a binary operation. As a generalization of this notion, and in order to encompass wider forms of binary interactions among the elements of a set, the following definition was proposed in [24].

**Definition 1.** *Let E be a non-empty set. A mapping from E×E into the power set P(E) of E, is called a law of synthesis on E. If the values of the law of synthesis are always singletons, then it is called a composition on E; otherwise, it is called a hypercomposition on E. A set equipped with a law of synthesis is called a magma.*

If " $\cdot$ " is a composition on a set $E$ and $A$, $B$ are subsets of $E$, then $A \cdot B$ signifies the set $\left\{ a \cdot b \mid (a,b) \in A \times B \right\}$, while if " $\cdot$ " is a hypercomposition, then, $A \cdot B = \bigcup_{(a,b) \in A \times B} a \cdot b$. $Ab$ and $aB$ have the same meaning as $A\{b\}$ and $\{a\}B$ respectively. In general, the singleton $\{a\}$ is identified with its member $a$.

**Theorem 1.** *Let A, $B \subseteq E$. If $A = \varnothing$ or $B = \varnothing$, then $AB = \varnothing$.*

Two special cases of hypercomposition are defined as follows: A hypercomposition is called *total* if $a \cdot b = E$ for every $a,b \in E$, and is called *degenerate* if $a \cdot b = \varnothing$ for every $a,b \in E$.

**Definition 2.** *Let E be a non-empty set equipped with a hypercomposition. The resulting magma is called a hypergroupoid if for every pair of elements a, b it holds that $ab \neq \varnothing$. Otherwise, it is called a partial hypergroupoid.*

**Definition 3.** *A law of synthesis (x,y) → x·y on a non- empty set E is called associative if*

$$(x \cdot y) \cdot z = x \cdot (y \cdot z)$$

*for all x, y, z in E. A magma endowed with an associative law of synthesis is called an associative magma.*

If the law of synthesis is a composition, the associative magma is referred to as a *semigroup*. If, however, the law of synthesis is a hypercomposition and the structure is an associative hypergroupoid, it is called a *semihypergroup*.

**Definition 4.** *A law of synthesis (x,y) → x·y on a non-empty set E is called reproductive if for all x in E, it holds that*

$$E \cdot x = x \cdot E = E$$

*A magma with a reproductive law of synthesis is called a reproductive magma.*

If the law of synthesis is a composition, the reproductive magma is called a *quasigroup*. If, however, the law of synthesis is a hypercomposition and the structure is a reproductive hypergroupoid, it is referred to as a *quasihypergroup*.

The reproductive axiom in the case of an associative magma whose law of synthesis is a composition (i.e., a semigroup), is equivalent to the classical group axioms, namely the existence of the identity and the inverse element (see also [21,24]). In particular, the following Theorem holds:

**Theorem 2.** *In a semigroup G, the reproductive axiom is equivalent to the postulates:*

   i. *There exists an element* $e \in G$ *such that* $ea = a = ae$, *for all* $a \in G$.

   ii. *For each element* $a \in G$ *there exists an element* $a' \in G$ *such that* $a'a = e = aa'$.

**Proof.** First, we assume that the reproductive axiom holds in $G$.

*i. Existence of the Neutral (Identity) Element.*

Let $x \in G$. By the reproductive axiom, $x \in xG$. Consequently, there exists $e \in G$ such that $xe = x$. Next, let $y$ be an arbitrary element of $G$. Since the composition is reproductive $y \in Gx$; therefore, there exists $z \in G$ such that $y = zx$. Using associativity we have:

$$ye = (zx)e = z(xe)$$

Since we established $xe = x$, we get $z(xe) = zx = y$. Thus $ye = y$.

In an analogous way, there exists an element $\hat{e}$ such that $\hat{e}y = y$ for all $y \in G$. Then the equality $e = \hat{e}e = \hat{e}$ is valid. Therefore, there exists $e \in G$ such that $ea = a = ae$, for all $a \in G$.

*ii. Existence of the Symmetric (Inverse) Element.*

By hypothesis, the reproductive axiom holds in $G$, and we have already established the existence of the identity element, $e$, from part (i). Let $a \in G$. By the reproductive axiom $e \in aG$. Thus, there exists $a' \in G$, such that $e = aa'$. Also, by the reproductive axiom, $e \in Ga$. Therefore, there exists $a'' \in G$, such that $e = a''a$. We show that $a'$ and $a''$ coincide. Using the associativity of and the identity property ($ea' = a'$ and $a''e = a''$), we obtain the following chain of equalities:

$$a' = ea' = (a''a)a' = a''(aa') = a''e = a''.$$

Hence, $a'$ and $a''$ coincide. Therefore $aa' = e = a'a$.

Conversely, let us assume that the two postulates (i) and (ii) are fulfilled. It must be proved that the reproductive axiom is valid for $G$. That is, we must show that for every $a \in G$, $aG = G$ and $Ga = G$. Proof of $aG = G$. Since $aG \subseteq G$, for all $a \in G$, it only remains to be proved that $G \subseteq aG$. Suppose that $x \in G$. Then, using the identity element $e$ and the inverse element $a'$ (from postulates i and ii, respectively), we have: $x = ex = (aa')x$. Applying the associativity of $G$, $x = a(a'x)$. The product $a'x$ is an element of $G$. Thus, $x$ can be expressed as $a$ multiplied by some element of $G$, which means $x \in aG$. Hence $G \subseteq aG$. Therefore $aG = G$. Similarly, we prove $Ga = G$. □

Hence, by Theorem 2, it becomes evident that a group may be defined solely in terms of associativity and reproductivity. Thus, the following definition can be stated:

**Definition 5** ([24]). *An associative and reproductive magma is called a hypergroup if the law of synthesis is a hypercomposition, while it is called a group if the law of synthesis is a composition.*

One can find three equivalent definitions for the group and two for the hypergroup in [24].

**Theorem 3.** *In a hypergroup H, the result of the synthesis of any two elements is always non-empty.*

**Proof.** Suppose that $xy = \varnothing$ for some $x,y \in H$. By reproductivity:

$$H = xH = x(yH)$$

By associativity:

$$x(yH) = (xy)H = \varnothing H$$

By Theorem 1, $\varnothing H = \varnothing$. Hence $H = \varnothing$, a contradiction. □

From the preceding theorem, it follows that the set of axioms defining a hypergroup can be reduced, since the first axiom, i.e., the non-emptiness of the result of the hypercomposition, is implied by the remaining two, namely, the associativity and the reproductivity of the structure.

**Corollary 4.** *The result of the hypercomposition in any structure built upon hypergroups is always non-empty.*

The structures referred to in the above corollary include hyperrings [13], hyperringoids [14], hyperfields [13], and, of course, hypermodules [15], hypermoduloids [16], vector hyperspaces [17], etc., which are all structures defined over hypergroups.

In the aforementioned hypercompositional structures, the additive part is a hypergroup, whereas the multiplicative part is a semigroup. However, in 1982, R. Rota introduced a hypercompositional structure called the multiplicative hyperring, in which, contrary to the above cases, the additive part is an abelian group, while the multiplicative part is a semihypergroup [25].

In particular, the definition of this structure is the following:

**Definition 6.** *A triplet $\left(S,+,\cdot\right)$ is called a multiplicative hyperring if*

1. *$\left(S,+\right)$ is an abelian group*

2. *$\left(S,\cdot\right)$ is a semihypergroup*

3. *$a \cdot \left(b + c\right) \subseteq a \cdot b + a \cdot c, \quad \left(b + c\right) \cdot a \subseteq b \cdot a + c \cdot a, \text{ for all } a,b,c \in S$*

4. *$a \cdot \left(-b\right) = \left(-a\right) \cdot b = -\left(a \cdot b\right), \text{ for all } a,b \in S$*

This structure has been extensively investigated both by R. Rota [26–28] and by numerous other researchers, e.g. [29–40]. In what follows, we demonstrate that, as with hyperrings, the hypercomposition in multiplicative hyperrings also consistently yields a non-empty result.

More precisely, we will show that the second Axiom in the aforementioned definition is stronger than strictly required. Specifically, while the definition of a semihypergroup requires the product of any two elements to be a non-empty set, we prove that in the case of multiplicative hyperrings, this condition is not independent; rather, it is a consequence of the remaining axioms. Consequently, we demonstrate that the definition of a multiplicative hyperring can be simplified as follows:

**Definition 7.** *A triplet $(S, +, \cdot)$ is called a multiplicative hyperring if*

1. *$(S, +)$ is an abelian group*

2. *$(S, \cdot)$ is a non-degenerate associative magma*

3. *$a \cdot (b + c) \subseteq a \cdot b + a \cdot c$, $(b + c) \cdot a \subseteq b \cdot a + c \cdot a$, for all $a, b, c \in S$*

4. *$a \cdot (-b) = (-a) \cdot b = -(a \cdot b)$, for all $a, b \in S$*

To justify the redundancy of the non-emptiness condition, arising from the consideration of the revised Axiom 2, we examine the properties of a structure $S$ that satisfies the axioms of Definition 7. We begin by establishing the following auxiliary result:

**Lemma 5.** *If $S$ contains an element $w$ such that $wz = \varnothing$, for some $z$ in $S$, then $wS = \varnothing$.*

**Proof.** As $(S, +)$ is a group, the reproductive axiom ensures that $z + S = S$, for every $z$ in $S$. Therefore, we have:

$$wS = w(z + S) \subseteq wz + wS$$

Given that $wz = \varnothing$, and by applying Theorem 1, the right side of the inclusion becomes $\varnothing + wS = \varnothing$. Thus, $wS = \varnothing$. □

**Theorem 6.** *In $S$ the product of any two elements is non-empty.*

**Proof.** Assume $wv = \varnothing$ for some $w$, $v$ in $S$ holds. Then, by Lemma 5, it follows that $wS = \varnothing$. Next let $x$ be any element in $S$. Since the additive part of $S$ is a group, the reproductive axiom guarantees that $S = w + S$. Thus $x = w + z$, for some $z \in S$. Let $y$ be an element in $S$. Using the distributivity law, we can write:

$$xy \subseteq xS = (w + z)S \subseteq wS + zS$$

Since $wS = \varnothing$, applying Theorem 1, the above reduces to:

$$xy \subseteq \varnothing + zS = \varnothing$$

Since $xy \subseteq \varnothing$, it follows that $xy = \varnothing$. This is a contradiction, as the hypercomposition is non-degenerate. □

Consequently, Theorem 6 confirms that the simplified axiomatic framework presented in Definition 7 is sufficient and internally consistent, as it guarantees the non-emptiness of the product over the entire set $S$ without the need for an explicit non-emptiness postulate.

Next, we proceed to demonstrate that the result of the hyperoperation is consistently non-empty in various other hypercompositional structures as well.

**Definition 8.** *A law of synthesis $(x, y) \rightarrow x \cdot y$ on a non-empty set $E$ is called weakly associative if the property,*

$$(x \cdot y) \cdot z \cap x \cdot (y \cdot z) \neq \varnothing$$

*is valid, for all $x$, $y$, $z$ in $E$. A magma whose law of synthesis is weakly associative is called a weakly associative magma.*

**Theorem 7.** *In a weakly associative magma, the hypercomposition yields a non-empty set.*

**Proof.** Suppose that $xy = \varnothing$, for some $x,y$ in $E$. Then,

$$(xy)z = \varnothing z = \varnothing , \text{ for any } z \in E.$$

Therefore,

$$(xy)z \cap x(yz) = \varnothing,$$

which is absurd. Hence, $xy$ is non-empty. □

Consequently, the following holds for H$v$-groups [12]:

**Corollary 8.** *The result of the law of synthesis of any two elements in H$v$-groups as well as in structures built upon them is always non-empty.*

**Definition 9.** *A law of synthesis* $(x,y) \rightarrow x \cdot y$ *on a non-void set E is called left inverted associative if the equality*

$$(x \cdot y) \cdot z = (z \cdot y) \cdot x$$

*is valid, for all elements x, y, z in E. It is called right inverted associative if*

$$x \cdot (y \cdot z) = z \cdot (y \cdot x)$$

*is valid, for all elements x, y, z in E. A magma which is equipped with a left inverted associative law of synthesis is called a left inverted associative magma, while one equipped with a right inverted associative law of synthesis is called a right inverted associative magma.*

**Definition 10.** [11] *A reproductive magma that satisfies the axiom of left inverted associativity is called a left almost-hypergroup (LA-hypergroup) if the law of synthesis is a hypercomposition and a left almost-group (LA-group) if the law of synthesis on the magma is a composition. Analogously, if the structure satisfies the axiom of right inverted associativity, it is termed a right almost-hypergroup (RA-hypergroup) or a right almost-group (RA-group), respectively.*

**Theorem 9.** *In every left almost-hypergroup (respectively, right almost-hypergroup) H the hypercomposition of any two elements of H is non-empty.*

**Proof**. Suppose that $H$ is a left almost-hypergroup and $xy = \varnothing$, for some $x,y$ in $H$. By reproductive axiom, we have $H = Hy$ and $H = Hx$. Then, applying the reproductivity and the left inverted associativity, we derive:

$$H = Hx = (Hy)x = (xy)H$$

Substituting the initial assumption $xy = \varnothing$ and applying Theorem 1, we get:

$$(xy)H = \varnothing H = \varnothing$$

Thus $H = \varnothing$, which is a contradiction since $H$ is non-empty.

Next assume that $H$ is a right almost-hypergroup and $xy = \varnothing$, for some $x,y$ in $H$ . Then by reproductivity, left inverted associativity and Theorem 1, we get:

$$H = yH = y(xH) = H(xy) = H\varnothing = \varnothing$$

which is again a contradiction. □

**Corollary 10.** *The result of the law of synthesis of any two elements in left (right) almost-hyperrings, left (right) almost-hypermodules etc. is always non-void.*

Every law of synthesis in a magma $E$ induces two new laws of synthesis [24]. If the law of synthesis is written multiplicatively, then the two induced laws are:

$$x/y = \{\, z \in E \mid x \in zy \,\} \quad \text{and} \quad y \backslash x = \{\, z \in E \mid x \in yz \,\}$$

which are called right division and left division, respectively. If the magma is commutative, then the right and left divisions coincide.

**Example 1.** *On the set $\mathbb{R}$ of real numbers, multiplication is a commutative law of synthesis. The inverse law is as follows*:

$$a\,/\,b = \begin{cases} \dfrac{a}{b} & \text{if } b \neq 0 \\[2mm] \varnothing & \text{if } b = 0,\, a \neq 0 \\[2mm] \mathbb{R} & \text{if } b = 0,\, a = 0 \end{cases}$$

*Here, the law of synthesis is a composition and the inverse law is a hypercomposition.*

**Theorem 11.** *In a reproductive magma $E$, the sets $x/y$ and $y\backslash x$ are non-empty for all $x,y \in E$. Conversely, if the sets $x/y$ and $y\backslash x$ are non-empty for every pair of elements $x,y$ in the magma, then the magma is reproductive.*

**Proof.** If $E$ is reproductive, then $Ey = E$ for all $y \in E$. Thus, for every $x \in E$ there exists $z \in E$, such that $x \in zy$. Therefore, $z \in x\,/\,y$ and thus $x/y \neq \varnothing$. Analogously, it follows that $y\backslash x$ is non-empty. Conversely, if $x/y \neq \varnothing$, for all $x,y \in E$, then, there exists $z \in E$, such that $x \in zy$, whence, $E \subseteq Ey$. Since $Ey \subseteq E$, for all $y \in E$, it follows that $Ey = E$. Likewise, $y\backslash x \neq \varnothing$, implies $yE = E$ for all $y \in E$. Hence, the magma is reproductive. □

**Corollary 12.** *In hypergroups, left almost-hypergroups, right almost-hypergroups, and $H_v$-groups, each of the two induced hypercompositions of any two elements is always non-empty.*

## 3. Reduction of Axioms in Polysymmetrical Hypergroups and Polysymmetrical Hyperrings

The notion of the *polysymmetrical hypergroup* was introduced into hypercompositional algebra in 1970 by J. Mittas [22], who was motivated by the theory of algebraically closed fields. Subsequently, this idea was extended, leading to the definition of many other polysymmetrical hypergroups [41–47]. Thus, a large family of such structures within hypergroups was formed, many of which emerged from the study of formal languages and automata using hypercompositional algebra [18-20, 48–50]. In honor of J. Mittas, C. Yatras named the polysymmetrical hypergroup originally defined in [22] the *M-polysymmetrical hypergroup*. The M-polysymmetrical hypergroup and the *M-polysymmetrical hyperring* were the subject of study by Yatras himself [51–55] as well as other researchers, e.g. [56–60]. The M-polysymmetrical hypergroup is a commutative hypergroup and serves as the prototype for the broader class of polysymmetrical hypergroups. In what follows, we present its formal definition.

**Definition 11.** *An M-polysymmetrical hypergroup H is a magma endowed with a hypercomposition that satisfies the following axioms:*

1. *Associativity:*
$$x(yz) = (xy)z, \text{ for all } x, y, z \in H$$

2. *Commutativity:*
$$xy = yx, \text{ for all } x, y \in H$$

3. *Existence of a neutral (or identity) element:*
   *There exists an element $e \in H$, such that*
$$x \in ex = xe, \text{ for all } x \in H$$

4. *Polysymmetry:*
   *For every $x \in H$ there exists at least one element $x' \in H$ such that*
$$xx' = x'x = e.$$
   *Such an element $x'$ is called symmetric of x and the set*
$$S(x) = \{x' \in H \mid xx' = x'x = e\}$$
   *is called the symmetric set of x (with respect to e).*

5. *Reversibility:*
$$\text{If } z \in xy \text{ and if } x' \in S(x), y' \in S(y), z' \in S(z) \text{ then } z' \in y'x'$$

By removing the commutativity axiom from the above definition, we obtain a more general structure, namely the *quasi-M-polysymmetrical hypergroup* [51], hereafter abbreviated as the *qMp-hypergroup*. Clearly, all results that hold for quasi-M-polysymmetrical hypergroups also hold for M-polysymmetrical hypergroups.

In the above definition, the axiom of associativity is explicitly stated, whereas the axiom of reproductivity does not appear. However, as the following theorem shows, this axiom follows from the remaining axioms of the definition.

**Theorem 13.** *The axiom of reproductivity holds in quasi-M-polysymmetrical hypergroups.*

**Proof.** Let $x \in H$. To establish reproductivity, we must show that $xH = H$ and $Hx = H$. Consider an arbitrary element $y \in H$. By the axiom of polysymmetry, there exists an element $x' \in S(x)$ such that $e \in xx'$. Using the property of the identity element and the associativity of the hyperoperation, we have:

$$y \in ey = (xx')y = x(x'y).$$

Thus, for $z \in x'y \subseteq H$, we have $y \in xz$. Consequently, every $y \in H$ is contained in the product of $x$ with some element of $H$, proving that $H \subseteq xH$. Since the reverse inclusion $xH \subseteq H$ holds by definition, we conclude that $xH = H$. An analogous argument establishes that $Hx = H$. Thus, the axiom of reproductivity is satisfied in $H$. □

The above Theorem rigorously demonstrates that the structure defined by Definition 11 is indeed a hypergroup.

**Proposition 14.** *In qMp-hypergroups, it holds that* $S(e) = \{e\}$.

**Proof.** Let $e' \in S(e)$. Then, by Axiom (4), $e'e = e$. Moreover, by Axiom (3), $e' \in e'e$. Hence, $e' = e$. Therefore, $S(e) = \{e\}$. □

**Corollary 15.** *In qMp-hypergroups, it holds that* $ee = e$.

At this point, let us recall some basic notions regarding identities in a hypergroup $H$. An element $e \in H$ is called a *right identity* if $x \in xe$ for all $x \in H$, and a *left identity* if $x \in ex$ for all $x \in H$. An element that is both a right and a left identity is simply called an *identity*. An identity is called *scalar* if $x = xe = ex$ for all $x \in H$, and *strong* if $xe = ex \subseteq \{e, x\}$ for all $x \in H$ [24]. An element $x \in H$, $x \neq e$ is called *attractive* if $e \in ex$ and $e \in xe$ [61]. The existence of such elements was revealed in the study of formal language and automata through hypercompositional algebra [18-20, 48-50]. The term attractive was chosen because such an element "attracts" the neutral element $e$ into the result of its hypercomposition with $e$.

**Proposition 16.** *The qMp-hypergroups do not have attractive elements.*

**Proof.** Let $H$ be a *qMp*-hypergroup and $x \in H$. We shall show that if $e \in ex$, then $x = e$. Pick $x' \in S(x)$. Then

$$e \in ex \Rightarrow ex' \subseteq (ex)x' = e(xx') = ee = \{e\}.$$

Hence $ex' = \{e\}$. Therefore

$$ex = (ex')x = e(x'x) = ee = \{e\},$$

so $ex = \{e\}$. Using Axiom (3), we also have $x \in ex$, whence $x = e$. □

**Proposition 17.** *In a qMp-hypergroup if* $y \in xy$, *then* $x = e$.

**Proof.** Let $H$ be an *qMp*-hypergroup with neutral element $e$ and let $x, y \in H$ such that $y \in xy$. Then we obtain the following chain of implications:

$$y \in xy \Rightarrow yy' \subseteq (xy)y' = x(yy') \Rightarrow e \in xe.$$

Therefore, by Proposition 16, it follows that $x = e$. □

**Proposition 18.** *In a qMp-hypergroup if* $x', x'' \in S(x)$, *then* $ex' = ex''$.

**Proof.** Let $t \in ex'$. Then $t \in (x''x)x' = x''(xx') = x''e$. Therefore $ex' \subseteq ex''$. Similarly, $ex'' \subseteq ex'$ and the Proposition follows. □

**Proposition 19.** *In a qMp-hypergroup H, given* $x \in H$ *and* $x' \in S(x)$, *we have*

$$S(x) = x'e.$$

**Proof.** Let $x' \in S(x)$. By the associativity of the hyperoperation and the property of the identity e, we have:

$$x(x'e) = (xx')e = ee = e$$

Thus $x'e \subseteq S(x)$. Conversely, suppose that $t \in S(x)$, then:

$$tx = e \Rightarrow (tx)x' = ex' \Rightarrow t(xx') = ex' \Rightarrow te = ex'$$

Under the assumption of the identity property $t \in te$, the equality $te = ex'$ imply $t \in ex'$. So $S(x) \subseteq ex'$. Combining the two inclusions, we obtain $S(x) = x'e$. □

**Proposition 20.** *In a qMp-hypergroup, if* $xe \cap ye \neq \varnothing$, *then*

    *i.*   $x, y \in S(x')$, *where* $x' \in S(x)$.

    *ii.*   $xe = ye$

**Proof.** *i.* Let $x' \in S(x)$. Then:

$$xe \cap ye \neq \varnothing \Rightarrow (x'x)e \cap (x'y)e \neq \varnothing \Rightarrow ee \cap (x'y)e \neq \varnothing \Rightarrow e \in (x'y)e$$

However, by Proposition 16, the elements of a qMp-hypergroup are not attractive. Hence, $x'y = e$ and so $y \in S(x')$, as required.

*ii.* If $xe \cap ye \neq \varnothing$, then, because of (i) $x, y \in S(x')$ and therefore, because of Proposition 18, $xe = ye$. □

Therefore, the Theorem holds:

**Theorem 21.** *In qMp-hypergroup H the family* $\{ xe \mid x \in H \}$ *is a partition of H.*

If $x \in H$, then the set $xe$, namely the class of the element $x$, will be denoted by $C(x)$.

**Proposition 22.** *In a qMp-hypergroup, if* $z \in xy$, *then* $C(z) = C(w)$ *for all* $w \in xy$.

**Proof.** Let $z \in xy$. Then we have the following implications:

$$z \in xy \Rightarrow zy' \subseteq (xy)y' \Rightarrow zy' \subseteq x(yy') \Rightarrow zy' \subseteq xe \Rightarrow zy' \subseteq ex \Rightarrow$$

$$\Rightarrow (zy')x' \subseteq (ex)x' \Rightarrow z(y'x') \subseteq e(xx') \Rightarrow z(y'x') \subseteq ee$$

By Corollary 15, $ee = e$, which leads to $z(y'x') = e$. Subsequently, the following implications hold:

$$(zy')x' = e \Rightarrow [(zy')x']x = ex \Rightarrow (zy')e = ex \Rightarrow e(zy') = ex \Rightarrow$$

$$\Rightarrow [e(zy')]y = (ex)y \Rightarrow ez = e(xy) \Rightarrow ez = \bigcup_{w \in xy} ew$$

Hence, by applying Theorem 21, we conclude the proposition. □

**Corollary 23.** *In a qMp-hypergroup, if* $xy \cap zw \neq \varnothing$, *then* $xy = zw$.

All the above theorems and propositions have been established without invoking the axiom of reversibility. We may therefore proceed to state and prove the following theorem:

**Theorem 24.** *In a quasi-M-polysymmetrical hypergroup, the axiom of reversibility is a consequence of the remaining axioms.*

**Proof.** Let $z \in H$. To establish reversibility, we show that if $z \in xy$ then for any $x' \in S(x)$, $y' \in S(y)$ and $z' \in S(z)$ the relation $z' \in y'x'$ holds. Given $z \in xy$, we derive the following implications:

$$z \in xy \Rightarrow zy' \subseteq x(yy') \Rightarrow zy' \subseteq xe \Rightarrow zy' \subseteq ex \Rightarrow z(y'x') \subseteq e(xx')$$

Since $xx' = e$ and $ee = e$ (by Corollary 15), it follows that $z(y'x') \subseteq e$. Consequently $e = z(y'x')$. Therefore $S(z) \cap y'x' \neq \varnothing$. Hence, by utilizing Proposition 22, we obtain $S(z) = y'x'$ and thus the axiom of reversibility is verified. □

With the above theorem, we have demonstrated that the number of axioms required to define an M-polysymmetrical hypergroup can be reduced. This leads to a more streamlined and rigorous formulation. We therefore arrive at the following refined definition:

**Definition 12.** *An M-polysymmetrical hypergroup H is a magma endowed with a hypercomposition that satisfies the following axioms:*

1.  *Associativity:*

$$x(yz) = (xy)z, \text{ for all } x, y, z \in H$$

2.  *Commutativity:*

$$xy = yx, \text{ for all } x, y \in H$$

3.  *Existence of a neutral (or identity) element:*
    *There exists an element $e \in H$, such that*

$$x \in ex = xe, \text{ for all } x \in H$$

4.  *Polysymmetry:*
    *For every $x \in H$ there exists at least one element $x' \in H$ such that*

$$xx' = x'x = e.$$

    *Such an element x' is called symmetric of x and the set*

$$S(x) = \{x' \in H \mid xx' = x'x = e\}$$

    *is called the symmetric set of x (with respect to e).*

*By omitting the commutativity axiom, one obtains the quasi M-polysymmetrical hypergroup.*

The M-polysymmetrical hypergroup serves as the basis for the definition of the *M-polysymmetrical hyperring*. The latter was introduced by J. Mittas [22] and has since been studied by him and by other authors (see, e.g. [54,56,57]). After the reduction of the axioms of the M-polysymmetrical hypergroup, the following new concise definition of the M-polysymmetrical hyperring may be stated:

**Definition 13.** *An M-polysymmetrical hyperring is an algebraic structure $(H,+,\cdot)$ where H is a non-empty set, "·" is an internal composition on H, and "+" is a hypercomposition on H which satisfies the axioms*:

**I.** *Multiplicative axiom*

$H = H^* \cup \{0\}$ *where $(H^*, \cdot)$ is a multiplicative semigroup and 0 is a bilaterally absorbing element of H, i.e., $0 \cdot x = x \cdot 0 = 0$, for all $x \in H$.*

*II.  Additive axioms*

   *1.  Associativity:*

$$x + (y + z) = (x + y) + z, \text{ for all } x, y, z \in H$$

   *2.  Commutativity:*

$$x + y = y + x, \text{ for all } x, y \in H$$

   *3.  Neutral element*

$$x \in 0 + x = x + 0, \text{ for all } x \in H$$

   *4.  Polysymmetry:*

   *For every $x \in H$ there exists at least one element $x' \in H$ such that*

$$x + x' = x' + x = 0 .$$

   *Such an element x' is called symmetric of x and the set*

$$S(x) = \left\{ x' \in H \ \middle| \ x + x' = x' + x = 0 \right\}$$

   *is called the symmetric set of x.*

*III.  Distributive axiom*

$$x \cdot (y + z) = x \cdot y + x \cdot z, \quad (y + z) \cdot x = y \cdot x + z \cdot x, \text{ for all } x, y, z \in H$$

## 4. Reduction of the Axioms in Hyperfields, Hyperrings, Hypermodules and Vector Hyperspaces

The concept of a *hyperfield* was introduced by M. Krasner in 1956 in his attempt to approximate a complete valued field via a sequence of such fields [62]. In 1966, during his lectures at the National Technical University of Athens, he further defined the notion of the *hyperring* and subsequently constructed the class of quotient hyperfields and hyperrings [13], which generalize the hyperfields he had originally introduced in [62].

The question of whether non-quotient hyperfields and hyperrings exist has been crucial for establishing the independence of hyperfield and hyperring theory from classical field and ring theory. Consequently, Krasner explicitly posed this question in [13].

The existence of non-quotient hyperfields and hyperrings was established in [63,64]. Later, another class of non-quotient hyperrings and hyperfields was introduced in [65] while the results of [63] were generalized in [66]. Today a large number of finite non-quotient hyperfields can be found in [67]. It is noteworthy that the investigation into the existence of non-quotient hyperfields gave rise to numerous questions within field theory itself [66–70].

The theory of hyperfields and hyperrings has been extensively developed, leading to dozens of research papers. For a comprehensive overview of this activity, see, for example, [63–107]. Across the literature on hyperfields, the structure is defined according to the axiomatization introduced by Krasner in [13]. We next present Krasner's definition as it appears in [13]:

**Definition 14.** *A hyperfield is an algebraic structure (H,+,·) where H is a non-empty set, "·" is an internal composition on H, and "+" is a hypercomposition on H which satisfies the axioms*:

*I.  Multiplicative axiom*

   $H = H^* \cup \{0\}$ *where $(H^*, \cdot)$ is a multiplicative group and 0 is a bilaterally absorbing element of*

$H$, i.e., $0 \cdot x = x \cdot 0 = 0$, for all $a \in H$.

**II. Additive axioms**

    1. associativity:

$$x + (y + z) = (x + y) + z, \text{ for all } x, y, z \in H$$

    2. commutativity:

$$x + y = y + x, \text{ for all } x, y \in H$$

    3. for every $x \in H$ there exists one and only one $x' \in H$ such that $0 \in x + x'$.
        $x'$ is written $-x$ and called the opposite of $x$; moreover, instead of $x + (-y)$ we write $x - y$.

    4. reversibility:

$$\text{if } z \in x + y, \text{ then } x \in z - y$$

**III. Distributive axiom**

$$z \cdot (x + y) = z \cdot x + z \cdot y, \quad (x + y) \cdot z = x \cdot z + y \cdot z, \text{ for all } x, y, z \in H$$

If the multiplicative axiom **I** is replaced by the axiom:

**I′.** $H = H^* \cup \{0\}$ where $(H^*, \cdot)$ is a multiplicative semigroup and $0$ is a bilaterally absorbing element of $H$

then, a more general structure is obtained which is called *hyperring* [13].

**Theorem 25.** *A non-empty set $H$ enriched with the additive axioms II is a hypergroup.*

**Proof.** Associativity holds by hypothesis. Next, assume that there exist elements $x, y \in H$ such that $x + y = \varnothing$. This leads to the following chain of implications:

$$x + y = \varnothing \Rightarrow (x + y) - y = \varnothing - y \Rightarrow x + (y - y) = \varnothing \Rightarrow x + 0 = \varnothing \Rightarrow x = \varnothing$$

This, however, contradicts the fact that $H$ is a non-empty set. Consequently, the result of the hypercomposition of any two elements in $H$ is always non-empty. Therefore, for every $y$, $w$ in $H$, there exists $x$ in $H$ such that $x \in w - y$ or, by reversibility, $w \in y + x$. This implies that $H \subseteq y + H$. Since the inclusion $y + H \subseteq H$ is given by the definition of the hypercomposition, we obtain $y + H = H$. Thus, the reproductive law is satisfied, and $H$ is a hypergroup. □

    The additive hypergroup of the hyperfield was termed a *canonical hypergroup* by J. Mittas as early as 1970 [108]; Mittas subsequently devoted a significant portion of his research to its systematic investigation (see, e.g. [108–111]). It has been established that canonical hypergroups are intimately linked to join hypergroups; specifically, a canonical hypergroup is essentially a join hypergroup characterized by the presence of a scalar neutral element [112]. Furthermore, a non-commutative canonical hypergroup is referred to as a *quasicanonical hypergroup* [113, 114] or, equivalently, a *polygroup* [115, 116].

**Theorem 26.** [63] *In a canonical hypergroup $H$ the identity $x + 0 = x$ holds for all $x$ in $H$.*

**Proof.** By definition, $0 \in x - x$; therefore, by reversibility, $x \in x + 0$. Next, let $w \in H$ be any element such that $w \in x + 0$. It follows from reversibility that $0 \in w - x$. According to Axiom

3, the opposite element is unique; consequently, $x=w$. Thus, $x+0=\{x\}$, which completes the proof. □

**Theorem 27.** *In a canonical hypergroup H the axiom of reversibility is equivalent to the validity of the equality:*

$$-(z+w)=-z-w$$

*for all z,w in H.*

**Proof.** Assume that the axiom of reversibility holds and consider an element $x \in -(z+w)$. By virtue of associativity and the reversibility property, the following chain of implications can be derived:

$$0 \in x+(z+w) \Rightarrow 0 \in (x+z)+w \Rightarrow -w \in x+z \Rightarrow x \in -z-w,$$

thereby proving that $-(z+w) \subseteq -z-w$. The reverse inclusion $-z-w \subseteq -(z+w)$ is established through a similar sequence of implications:

$$x \in -z-w \Rightarrow -z \in x+w \Rightarrow 0 \in (x+w)+z \Rightarrow 0 \in x+(w+z) \Rightarrow x \in -(z+w)$$

Conversely, assume the identity $-(z+w)=-z-w$ holds for all $z,w$ in $H$. We shall verify the reversibility condition: $z \in x+y \Rightarrow x \in z-y$. The premise $z \in x+y$ is equivalent to $0 \in z-(x+y)$ which, by applying the equality $-(x+y)=-x-y$, becomes $0 \in z+(-x-y)$ or $0 \in (x-y)-z$. The relation $0 \in (x-y)-z$ directly yields $z \in x-y$, which confirms the reversibility. □

The above theorem establishes a fundamental equivalence for canonical hypergroups, an equivalence whose significance becomes even more pronounced in the context of hyperfields. Since the identity $-(x+y)=-x-y$ holds in every hyperfield as a direct consequence of the distributivity axiom, the following can be stated:

**Theorem 28.** *In a hyperfield H, reversibility is a consequence of the remaining axioms.*

Theorem 28 therefore enables an axiomatic simplification; specifically, the reversibility condition, traditionally included in the definition of the additive structure of the hyperfield, becomes redundant. Consequently, this allows for the statement of a more concise definition of a hyperfield:

**Definition 15.** *A hyperfield is an algebraic structure (H,+,·) where H is a non-empty set, "·" is an internal composition on H, and "+" is a hypercomposition on H which satisfies the axioms:*

**I.** *Multiplicative axiom*

$H = H^* \cup \{0\}$ *where* $(H^*,\cdot)$ *is a multiplicative group and 0 is a bilaterally absorbing element of H, i.e.,* $0 \cdot x = x \cdot 0 = 0$, *for all* $a \in H$.

**II.** *Additive axioms*

*1. associativity:*

$$x+(y+z)=(x+y)+z, \text{ for all } x,y,z \in H$$

2. *commutativity:*
$$x + y = y + x, \text{ for all } x, y \in H$$

3. *for every $x \in H$ there exists one and only one $x' \in H$ such that $0 \in x + x'$.*
   *$x'$ is written $-x$ and called the opposite of $x$; moreover, instead of $x + (-y)$ we write $x - y$.*

**III.** *Distributive axiom*
$$z \cdot (x + y) = z \cdot x + z \cdot y, \quad (x + y) \cdot z = x \cdot z + y \cdot z, \text{ for all } x, y, z \in H$$

The definition of a hyperring is formulated in an analogous way provided that a multiplicative identity is present.

A fundamental question arises when a hyperring lacks a multiplicative identity: Can any hyperring be embedded into a hyperring that possesses an identity?

For ordinary rings, a standard construction achieving this is the Dorroh extension. However, this construction cannot be applied when the additive law is a hypercomposition rather than a composition.

Let us examine why this construction fails. Assume that $R$ is a hyperring and let $U$ be the product set $\mathbb{Z} \times R$. Two pairs $(n, x)$, $(m, y)$ are regarded equal if and only if $n = m$ and $x = y$. Define addition by:
$$(n, x) + (m, y) = \left\{ (n + m, z) \,\middle|\, z \in x + y \right\}$$

It is straightforward to verify that $U$, endowed with this hypercomposition, becomes a canonical hypergroup with neutral element $(0, 0)$, where the opposite of any element $(n, x)$ is $(-n, -x)$. However, proceeding with the Dorroh extension for multiplication leads to structural complications. The multiplication is defined by:
$$(n, x) \cdot (m, y) = \left\{ (nm, z) \,\middle|\, z \in ny + mx + xy \right\}$$

Since the term $ny + mx + xy$ represents a subset of $R$ rather than a single element, the multiplication is a hypercomposition and not a well-defined single-valued composition. Consequently, the Dorroh extension of a hyperring is not a hyperring.

J. Mittas, in the context of studying polynomials in one variable with coefficients in a hyperring, introduced a structure where both addition and multiplication are hypercompositions, terming this structure a *superring* [117, 118]. Specifically, a superring is a hypercompositional structure that generalizes the concept of a hyperring, as its multiplicative component forms a semihypergroup rather than a semigroup. Nevertheless, the Dorroh extension of a hyperring fails to be even a superring, as the multiplication defined above is not associative. Indeed, consider the product:
$$\left[ (n, x) \cdot (m, y) \right] \cdot (k, z) = \left\{ (nm, w) \,\middle|\, w \in ny + mx + xy \right\} \cdot (k, z) =$$
$$= \bigcup_{w \in ny + mx + xy} (nm, w) \cdot (k, z) = \bigcup_{w \in ny + mx + xy} \left\{ ((nm)k, v) \,\middle|\, v \in (nm)z + kw + wz \right\}$$

But

$$\bigcup_{w \in ny+mx+xy} \left\{ \big((nm)k, v\big) \,\Big|\, v \in (nm)z + kw + wz \right\} \subseteq$$

$$\subseteq \left\{ \big((nm)k, v\big) \,\Big|\, v \in (nm)z + k\big(ny+mx+xy\big) + \big(ny+mx+xy\big)z \right\}$$

because

$$\left\{ \big((nm)k, v\big) \,\Big|\, v \in (nm)z + k\big(ny+mx+xy\big) + \big(ny+mx+xy\big)z \right\} =$$

$$= \bigcup_{w,u \in ny+mx+xy} \left\{ \big((nm)k, v\big) \,\Big|\, v \in (nm)z + kw + uz \right\}$$

Therefore

$$\Big[(n,x)\cdot(m,y)\Big]\cdot(k,z) \subseteq$$

$$\subseteq \left\{ \big(n(mk), v\big) \,\Big|\, v \in (nm)z + (kn)y + (km)x + k\big(xy\big) + n\big(yz\big) + m\big(xz\big) + \big(xy\big)z \right\}$$

Similarly, for the right-side association:

$$(n,x)\cdot\Big[(m,y)\cdot(k,z)\Big] = (n,x)\cdot\left\{ (mk,w) \,\Big|\, w = mz+ky+yz \right\} =$$

$$= \bigcup_{w \in mz+ky+yz} (n,x)\cdot(mk,w) = \bigcup_{w \in mz+ky+yz} \left\{ \big(n(mk), v\big) \,\Big|\, v \in nw + (mk)x + xw \right\} \subseteq$$

$$\subseteq \left\{ \big(n(mk), v\big) \,\Big|\, v \in n\big(mz+ky+yz\big) + (mk)x + x\big(mz+ky+yz\big) \right\} =$$

$$= \left\{ \big(n(mk), v\big) \,\Big|\, v \in (nm)z + (kn)y + (km)x + k\big(xy\big) + n\big(yz\big) + m\big(xz\big) + \big(xy\big)z \right\}$$

Consequently, there is no guarantee that associativity—or even weak associativity—holds for this specific hypercomposition. It appears, therefore, that the corresponding theorem for rings, that any ring can be embedded into a ring with identity, does not hold for hyperrings via the standard construction methods.

The implications of Theorem 27 are not restricted to hyperfields and unitary hyperrings but extend to hypermodules and vector hyperspaces. These structures have been extensively studied by numerous researchers, yielding significant results (e.g. [119-135]) and various applications (e.g. [15,20,21]). The classical definition of a hypermodule is as follows [15]:

**Definition 16.** *A left hypermodule over a unitary hyperring $P$ is a canonical hypergroup $M$ with an external composition $\big(a,m\big) \rightarrow am$, from $P \times M$ to $M$ satisfying the conditions:*

***i.*** $\quad a\big(m+n\big) = am + an,$

***ii.*** $\quad \big(a+b\big)m = am + bm,$

***iii.*** $\quad \big(ab\big)m = a\big(bm\big),$

***iv.*** $\quad 1m = m$ and $\ 0m = 0$

*for all $a,b \in P$ and all $m,n \in M$.*

The *right hypermodule* is defined analogously. Furthermore, a hypermodule over a hyperfield is referred to as a *vector hyperspace*.

It is evident that, due to Axiom (i) and the fact that $P$ is a unitary hyperring, the equality

$$-(m+n) = -m-n,$$

holds. Consequently, the requirement that $M$ must be a canonical hypergroup in the definition of a hypermodule is redundant; it suffices for $M$ to be a commutative hypergroup possessing a neutral element 0 such that $x+0 = x$, for every $x \in M$ and for each $x \in M$ there exists a unique element $x' \in M$ such that $0 \in x+x'$. This observation similarly applies to vector hyperspaces. Such a hypergroup—not necessarily commutative—was termed a *normal hypergroup* by Marty [2, 3], who also referred to it as *completely regular* (*complètement régulier* [1]). A normal hypergroup further endowed with the axiom of reversibility constitutes a quasicanonical hypergroup (or polygroup). These considerations lead to the following refined definition:

**Definition 17.** *A left hypermodule over a unitary hyperring $P$ is a commutative normal hypergroup $(M,+)$ equipped with an external composition $(a,m) \rightarrow am$ from $P \times M$ into $M$ that satisfies condition (i) through (iv) of Definition 16.*

Whenever Definition 17 is invoked, it should be explicitly noted that the additive hypergroup involved is not merely normal: it is in fact canonical. For this reason, this definition must be accompanied by the theorem that follows.

**Theorem 29.** *The additive hypergroup of a hypermodule over a unitary hyperring is a canonical hypergroup.*

As proven in Theorem 16 of [20], the direct sum of hypermodules is not a hypermodule. This fact necessitated the introduction of the *weak hypermodule* ([20], Definition 10), which satisfies the inclusion condition $(a+b)m \subseteq am + bm$, in place of the equality in Axiom (ii). Within the framework of weak hypermodules—which arise naturally as the direct sum of hypermodules—it is evident that the additive structure of $M$ involved in Definition 17 is fully sufficient for the weak hypermodule to be well-defined.

## 5. Discussion

The clarification of the notion of hypercomposition and the reduction of axioms in the definitions of hypergroups, polysymmetrical hypergroups, hyperfields, unitary hyperrings, multiplicative hyperrings, hypermodules, vector hyperspaces, and various other hypercompositional structures serve not only to strengthen the logical foundations of these systems but also to enhance their applicability in algorithmic contexts. Simplifying the axiomatic framework fosters conceptual clarity and ensures mathematical rigor by explicitly distinguishing primitive, essential assumptions from derivable properties.

Furthermore, such axiomatic refinement is pivotal in the field of computational mathematics. By minimizing the set of required axioms, it becomes possible to design more streamlined algorithms characterized by reduced computational overhead. Consequently, the exhaustive construction and classification of finite hypercompositional structures can be achieved more effectively. This methodological approach was instrumental in the development of algorithms that enabled the complete construction of all hyperfields of order seven [67].

CORE DEPARTMENT, NATIONAL AND KAPODISTRIAN UNIVERSITY OF ATHENS
*Email address:* chrmas@uoa.gr *or* ch.massouros@gmail.com